\documentclass[reqno]{amsart}

\usepackage{amsmath,amssymb,dsfont}

\usepackage{fullpage}
\usepackage{tikz-cd}

\usepackage{hyperref} 

\usepackage[nobysame,alphabetic,initials,msc-links,lite]{amsrefs}

\DefineSimpleKey{bib}{how}
\renewcommand{\PrintDOI}[1]{%
  \href{http://dx.doi.org/#1}{{\tt DOI:#1}}%
}
\renewcommand{\eprint}[1]{#1}
\BibSpec{book}{%
    +{}  {\PrintPrimary}                {transition}
    +{.} { \PrintDate}                  {date}
    +{.} { \textit}                     {title}
    +{.} { }                            {part}
    +{:} { \textit}                     {subtitle}
    +{,} { \PrintEdition}               {edition}
    +{}  { \PrintEditorsB}              {editor}
    +{,} { \PrintTranslatorsC}          {translator}
    +{,} { \PrintContributions}         {contribution}
    +{,} { }                            {series}
    +{,} { \voltext}                    {volume}
    +{,} { }                            {publisher}
    +{,} { }                            {organization}
    +{,} { }                            {address}
    +{,} { }                            {status}
    +{,} { \PrintDOI}                   {doi}
    +{,} { \eprint}        {eprint}
    +{,} { \PrintISBNs}                 {isbn}
    +{}  { \parenthesize}               {language}
    +{}  { \PrintTranslation}           {translation}
    +{;} { \PrintReprint}               {reprint}
    +{.} { }                            {note}
    +{.} {}                             {transition}
    +{}  {\SentenceSpace \PrintReviews} {review}
}
\BibSpec{article}{%
    +{}  {\PrintAuthors}                {author}
    +{,} { \textit}                     {title}
    +{.} { }                            {part}
    +{:} { \textit}                     {subtitle}
    +{,} { \PrintContributions}         {contribution}
    +{.} { \PrintPartials}              {partial}
    +{,} { }                            {journal}
    +{}  { \textbf}                     {volume}
    +{}  { \PrintDatePV}                {date}
    +{,} { \issuetext}                  {number}
    +{,} { \eprintpages}                {pages}
    +{,} { }                            {status}
    +{,} { \PrintDOI}                   {doi}
    +{,} { \eprint}        {eprint}
    +{}  { \parenthesize}               {language}
    +{}  { \PrintTranslation}           {translation}
    +{;} { \PrintReprint}               {reprint}
    +{.} { }                            {note}
    +{.} {}                             {transition}
    +{}  {\SentenceSpace \PrintReviews} {review}
}
\BibSpec{collection.article}{%
    +{}  {\PrintAuthors}                {author}
    +{,} { \textit}                     {title}
    +{.} { }                            {part}
    +{:} { \textit}                     {subtitle}
    +{,} { \PrintContributions}         {contribution}
    +{,} { \PrintConference}            {conference}
    +{}  {\PrintBook}                   {book}
    +{,} { }                            {booktitle}
    +{,} { \PrintDateB}                 {date}
    +{,} { pp.~}                        {pages}
    +{,} { }                            {publisher}
    +{,} { }                            {organization}
    +{,} { }                            {address}
    +{,} { }                            {status}
    +{,} { \PrintDOI}                   {doi}
    +{,} { \eprint}        {eprint}
    +{}  { \parenthesize}               {language}
    +{}  { \PrintTranslation}           {translation}
    +{;} { \PrintReprint}               {reprint}
    +{.} { }                            {note}
    +{.} {}                             {transition}
    +{}  {\SentenceSpace \PrintReviews} {review}
}
\BibSpec{misc}{%
  +{}{\PrintAuthors}  {author}
  +{,}{ \textit}      {title}
  +{.}{ }             {how}
  +{}{ \parenthesize} {date}
  +{,} { available at \eprint}        {eprint}
    +{,} { \PrintDOI}                   {doi}
  +{,}{ available at \url}{url}
  +{,}{ }             {note}
  +{.}{}              {transition}
}

\numberwithin{equation}{section}

\newtheorem{theorem}{Theorem}[section]
\newtheorem{corollary}[theorem]{Corollary}
\newtheorem{lemma}[theorem]{Lemma}
\newtheorem{proposition}[theorem]{Proposition}

\theoremstyle{remark}
\newtheorem{remark}[theorem]{Remark}

\theoremstyle{definition}

\newcommand{\bp}{\begin{proof}}
\newcommand{\ep}{\end{proof}}

\newcommand{\C}{\mathbb{C}}
\newcommand{\N}{\mathbb{N}}
\newcommand{\T}{\mathbb{T}}
\newcommand{\Z}{\mathbb{Z}}

\newcommand{\A}{\mathcal{A}}
\newcommand{\CC}{\mathcal{C}}
\newcommand\DD{\mathcal D}
\newcommand{\BB}{\mathcal{B}}

\newcommand{\OO}{\mathcal{O}}

\newcommand{\U}{\mathcal{U}}
\newcommand{\Mult}{\mathcal{M}}

\newcommand{\TL}{\mathcal{TL}}

\newcommand{\un}{\mathds{1}}

\newcommand\QDG{\OO_c(\coD(G))}
\newcommand\QDS{\OO_c(\coD(\SU_q(2)))}
\mathchardef\mhyph="2D
\newcommand{\indcat}[1]{{\mathrm{ind}\mhyph#1}}
\newcommand{\indC}{\indcat{\CC}}
\newcommand{\Dcen}{\mathcal{Z}}
\newcommand{\ZC}{{\Dcen(\indC)}}

\newcommand{\Zreg}{Z_{\mathrm{reg}}}

\newcommand{\norm}[1]{\left\|#1\right\|}

\newcommand{\coD}{\hat{\mathcal D}}
\newcommand{\hlf}{{1/2}}
\newcommand\op{\mathrm{op}}

\DeclareMathOperator{\Dhat}{\hat\Delta}

\DeclareMathOperator{\Mor}{Mor}
\DeclareMathOperator{\End}{End}
\DeclareMathOperator{\Hom}{Hom}
\DeclareMathOperator{\Tr}{Tr}

\DeclareMathOperator{\Irr}{Irr}

\DeclareMathOperator{\Rep}{Rep}

\DeclareMathOperator{\onb}{onb}

\newcommand\SU{\mathrm{SU}}

\DeclareMathOperator{\modl}{mod}
\newcommand{\Qmod}{Q\mhyph\mathrm{mod}}
\newcommand{\modQ}{\mathrm{mod}\mhyph Q}
\newcommand{\QmodQ}{Q\mhyph\mathrm{mod}\mhyph Q}

\DeclareMathOperator{\Tub}{Tub}

\begin{document}

\title{A few remarks on the tube algebra of a monoidal category}

\author{Sergey Neshveyev}
\thanks{Supported by the European Research Council under the European Union's Seventh
Framework Programme (FP/2007-2013)/ ERC Grant Agreement no. 307663
}
\address{Department of Mathematics, University of Oslo, P.O. Box 1053 Blindern, NO-0316 Oslo, Norway}
\email{sergeyn@math.uio.no}

\author{Makoto Yamashita}
\thanks{Supported by JSPS KAKENHI Grant Number 25800058. Partially supported by the Danish National Research Foundation through the Centre for Symmetry and Deformation (DNRF92)}
\address{Department of Mathematics, Ochanomizu University, Otsuka
2-1-1, Bunkyo, 112-8610 Tokyo, Japan}
\email{yamashita.makoto@ocha.ac.jp}

\date{November 23, 2015; minor update March 10, 2018}

\keywords{C$^*$-tensor category; tube algebra; Morita equivalence}

\subjclass[2010]{Primary: 81R15; Secondary: 18D10, 46L89}

\begin{abstract}
We prove two results on the tube algebras of rigid C$^*$-tensor categories. The first is that the tube algebra of the representation category of a compact quantum group $G$ is a full corner of the Drinfeld double of $G$. As an application we obtain some information on the structure of the tube algebras of the Temperley--Lieb categories $\TL(d)$ for $d>2$. The second result is that the tube algebras of weakly Morita equivalent C$^*$-tensor categories are strongly Morita equivalent. The corresponding linking algebra is described as the tube algebra of the $2$-category defining the Morita context.
\end{abstract}

\maketitle

\section{Introduction}

The present work is motivated by the representation theory of the fusion algebras of C$^*$-tensor categories, a study of which was initiated by Popa and Vaes~\cite{MR3406647}, and by the authors~\cite{arXiv:1501.07390}. It was then realised by Ghosh and C.~Jones~\cite{MR3447719} that this theory can also be understood through the representations of the tube algebras, that contain the fusion algebras as nonfull corners. In this note our goal is to clarify the connections between the tube algebras and some of the constructions and problems studied in the above mentioned papers.

\smallskip

The Drinfeld double of a  Hopf algebra and its categorical counterpart, the Drinfeld centre of a monoidal category, are very powerful
general tools for producing modular categories, which are categories with nontrivial braiding symmetry. They played an important role in the
development of the theory of quantum groups during the 90's, as such
structures have rich connections with topological and conformal field
theories in mathematical physics.

Parallel to this development, in the framework of subfactor theory,
Ocneanu developed another formalism of a `quantum double' construction
based on the notion of asymptotic inclusion, see, for
example,~\cite{MR1316301}. He also introduced tube algebras
as a crucial tool to analyse certain systems of bimodules associated
with asymptotic inclusions, which led to a new construction of
$3$-dimensional topological quantum field theories.

The relation between the Drinfeld centre and the subfactor theoretic
construction was subsequently clarified through the effort of many people: Longo and Rehren~\cite{MR1332979}, Izumi~\cite{MR1782145},
Masuda~\cite{MR1442437}, and M\"{u}ger~\cite{MR1966525}, to name just a
few.  We should note that both notions were mostly
analysed within the framework of fusion categories, which imposes
finiteness assumption on the number of simple objects.  However, many constructions can be carried out without this
restriction. The
asymptotic inclusion is generalised to the notion of symmetric
enveloping (SE) inclusion due to Popa~\cite{MR1302385}, while the notions of Drinfeld centre and tube algebra do not require finiteness from the beginning. One difference, though, is that in the nonfusion case the unitary Drinfeld centre in the usual sense may be too small to be interesting~\citelist{\cite{MR1444286}\cite{arXiv:1501.07390}}, and we should rather consider the centres of the corresponding ind-categories. Taking this into account, the relation between the SE-inclusions and the Drinfeld centres was clarified in our previous work~\cite{arXiv:1501.07390}. As for the tube algebras, their representations categories are still
equivalent to the Drinfeld centres of the ind-categories~\cite{arXiv:1511.07329}. We should also note that in the
theory of planar algebras, an analogue of the tube algebra was formulated
as the annular category of a planar algebra~\cite{MR1929335}.

One new phenomenon in the non-fusion case is that the Drinfeld centre of an ind-category is no longer guaranteed to be semisimple, even if the
original category is. Consequently, the set of equivalence classes of simple objects
starts to have a much finer structure as a topological space with the
Fell topology, defined as the unitary dual of an appropriate algebra. This topological structure has direct applications to the
approximation properties of subfactors, see~\citelist{\cite{MR3406647}\cite{arXiv:1501.07390}\cite{MR3447719}\cite{MR3406863}\cite{MR3464395}}. Different algebras describing the same category can in principle define different topologies. One way
of dealing with this problem is to establish a strong Morita equivalence of such algebras,
which induces a natural homeomorphism with respect to the Fell
topology~\cite{MR0367670}. Our goal is to prove two results of this type.

\smallskip

The note is organised as follows. In Section~\ref{sec:prelim} we recall some basic notions and fix our conventions.

In Section~\ref{sec:tub-alg-rep-cat} we consider the representation category of a compact quantum group. As discussed above, the Drinfeld centre of the corresponding ind-category is described by two algebras, the Drinfeld double and the tube algebra. If we consider the
$q$-deformation of a compact semisimple Lie group, then the unitary dual of the Drinfeld double has a close resemblance to the unitary dual of the complexification of the original group~\citelist{\cite{MR1804864}\cite{arXiv:1410.6238}}. On the other
hand, the tube algebra seems to be more appropriate for combinatorial
analysis~\cite{MR3464395}. Through a detailed analysis of matrix coefficients of representations, we show that the tube algebra is contained in the Drinfeld double as a full corner. Although the result is not surprising,  even for finite groups an explicit connection between the two algebras seems to be missing in the literature. As an application we consider the Temperley--Lieb categories $\TL(d)$ and using the results of Pusz~\cite{MR1213303} we get some information on the structure of the tube algebra of $\TL(d)$ for $d>2$. We note that the representations of this algebra were studied by V.~Jones and Reznikoff~\citelist{\cite{MR1929335}\cite{MR2274519}} and Ghosh and C.~Jones~\cite{MR3447719}.

In Section~\ref{sec:mor-eqv-cat-tub-alg} we consider weakly Morita equivalent categories. By making use of the $2$-categorical formulation, we prove strong Morita equivalence of the corresponding tube algebras. Again, the result is known for fusion categories and is not surprising in general, since the corresponding Drinfeld centres are monoidally equivalent~\citelist{\cite{MR1822847}\cite{arXiv:1501.07390}}. But even for fusion categories an explicit imprimitivity bimodule is not available in the literature. We also briefly discuss an alternative approach based on `regular half-braidings with coefficients' inspired by M\"{u}ger's description of the tube algebra of a fusion category~\cite{MR1966525}. This complements our analysis of the fusion algebras of weakly Morita equivalent categories in~\cite{arXiv:1501.07390}, where we showed that certain approximation properties are preserved, yet the fusion algebras are not strongly Morita equivalent.

\smallskip
\noindent{\bf Acknowledgement.} Part of this work was carried out during the authors' participation in the Graduate School ``Topological Quantum Groups'' held at the Mathematical Research and Conference Center in B\c{e}dlewo. We would like to thank the organisers for their hospitality. We would also like to thank S.~Vaes for fruitful correspondence.

\bigskip

\section{Preliminaries}
\label{sec:prelim}

\subsection{Tensor categories}
We study essentially small strict rigid C$^*$-tensor categories, for which we follow the same conventions as in \cite{arXiv:1501.07390}. Therefore we assume that for any such category $\CC$ the tensor unit $\un$ is simple (which forces $\CC$ to be semisimple), and $\CC$ is closed under taking subobjects and finite direct sums. Let us recall a few basic notions and facts that we will repeatedly use.

\emph{Rigidity} means that every object $X$ in $\CC$ has a dual, that is,  there is an object $\bar{X}$ in $\CC$ and morphisms $R \in \CC(\un, \bar{X} \otimes X)$ and $\bar{R} \in \CC(\un, X \otimes \bar{X})$ satisfying the \emph{conjugate equations}
$$
  (\iota_{\bar{X}} \otimes \bar{R}^*) (R \otimes \iota_{\bar{X}})= \iota_{\bar{X}},\ \
  (\iota_X \otimes R^*) (\bar{R} \otimes \iota_X) = \iota_X.
$$
The \emph{intrinsic dimension} of $X$ is defined by
$$
d(X)=d^\CC(X) = \min_{(R, \bar{R})} \norm{R} \norm{\bar{R}},
$$
where $(R, \bar{R})$ runs over the solutions of the conjugate equations for $X$. A solution $(R, \bar{R})$ satisfying $\norm{R} = \norm{\bar{R}} = d(X)^{\hlf}$ is called \emph{standard}, and such solutions are unique up to transformations of the form $(R,\bar R)\mapsto((T \otimes \iota)R, (\iota \otimes T)\bar{R})$ for unitary morphisms~$T$.  We often denote a choice of standard solution for the conjugate equations for $X$ as $(R_X, \bar{R}_X)$. When $\{X_i\}_{i \in I}$ is a parameterised family of objects in $\CC$, we write $(R_i, \bar{R}_i)$ instead of $(R_{X_i}, \bar{R}_{X_i})$. Similarly, for many other constructions we use index $i$ instead of~$X_i$, so for example we write $d_i$ for $d(X_i)$. If the family is self-dual, we also write $\bar{\imath}$ for the index corresponding to the dual of $X_i$.

The \emph{categorical trace} is the trace on the endomorphism ring $\CC(X)=\CC(X,X)$ of $X$ defined by
$$
\Tr_X(T) = R^*_X(\iota \otimes T)R_X = \bar{R}^*_X(T \otimes \iota)\bar{R}_X.
$$
It is independent of the choice of standard solutions $(R_X, \bar{R}_X)$.

For $X, Y \in \CC$ and a choice of standard solutions $(R_X, \bar{R}_X)$ and $(R_Y, \bar{R}_Y)$, we can define a linear anti-multiplicative map $\CC(X, Y) \to \CC(\bar{Y}, \bar{X})$, denoted by $T \mapsto T^\vee$, which is characterised by $(T \otimes \iota) \bar{R}_X = (\iota \otimes T^\vee) \bar{R}_Y$. This map can be also characterized by $(\iota \otimes T) R_X = (T^\vee \otimes \iota) R_Y$ and satisfies $T^{\vee *} = T^{* \vee}$ for the standard solutions of the conjugate objects chosen as $(R_{\bar{X}}, \bar{R}_{\bar{X}}) = (\bar{R}_X, R_X)$.

\subsection{Fusion algebra}

There are several $*$-algebras associated with $\CC$. The better known, and easier to define, is the \emph{fusion algebra} $\C[\Irr(\CC)]$. As a space it is spanned by the isomorphism classes $[U]$ of objects in $\CC$ and satisfies the relations $[U\oplus V]=[U]+[V]$. The product and involution are defined by $[U]\, [V]=[U\otimes V]$ and $[U]^*=[\bar U]$. If we fix representatives $\{U_s\}_{s\in\Irr(\CC)}$ of the isomorphism classes of simple objects, then clearly $\{[U_s]\}_{s}$ is a basis in $\C[\Irr(\CC)]$.

In general $\C[\Irr(\CC)]$ does not admit a universal C$^*$-completion, but it is still possible to define a completion which plays the role of a full C$^*$-algebra. This can be done in several equivalent ways~\citelist{\cite{MR3406647}\cite{arXiv:1501.07390}\cite{MR3447719}}. Let us recall our approach in~\cite{arXiv:1501.07390}. Consider the C$^*$-tensor category $\indC$ of ind-objects of $\CC$. We refer the reader to~\cite{arXiv:1501.07390} for the precise definition, but informally $\indC$ is obtained from $\CC$ by allowing infinite direct sums of objects. We then consider the unitary Drinfeld centre $\ZC$ of $\indC$, so the objects of $\ZC$ are pairs $(Z,c)$, where $Z$ is an ind-object and $c$ is a unitary half-braiding on~$Z$, that is, a collection of unitary isomorphisms $c_X\colon X\otimes Z\to Z\otimes X$ which is natural in $X\in\indC$, such that for all objects~$X$ and~$Y$ in $\indC$ we have
\begin{equation*} \label{eq:halfbr}
c_{X\otimes Y}=(c_X\otimes\iota_Y)(\iota_X\otimes c_Y).
\end{equation*}
Every object $(Z,c)$ defines a representation $\pi_Z=\pi_{(Z,c)}$ of $\C[\Irr(\CC)]$ as follows. The underlying Hilbert space of the representation is
$H_Z=\Mor_{\indC}(\un, Z)$ with the scalar product such that $(\xi,\zeta)\iota=\zeta^*\xi$. Then for an object $X$ in~$\CC$ and a vector $\xi \in \Mor_\indC(\un, Z)$, we put $\pi_Z([X])\xi$ to be the composition
\begin{multline*}
(\iota_Z \otimes \bar{R}_X^*)(c_X \otimes \iota_{\bar{X}})(\iota_X \otimes \xi \otimes \iota_{\bar{X}})\bar{R}_X\colon
\un \xrightarrow{\bar R_X} X \otimes \bar{X} \xrightarrow{\iota\otimes\xi\otimes\iota} X \otimes Z \otimes \bar{X} \\
\xrightarrow{c_X \otimes \iota} Z \otimes X \otimes \bar{X} \xrightarrow{\iota_Z \otimes \bar{R}_X^*} Z.
\end{multline*}
From this definition $\norm{\pi_Z([X])}$ is bounded by $d(X)$, and we denote the completion of $\C[\Irr(\CC)]$ with respect to the norm $\sup\|\pi_Z(\cdot)\|$, where the supremum is taken over all $(Z,c)\in\ZC$, by $C^*(\CC)$.

Denote by $x\mapsto x^\vee$ the $*$-anti-automorphism of $\C[\Irr(\CC)]$ characterised by $[U]^\vee = [\bar U]$. Then for every representation $\pi\colon \C[\Irr(\CC)]\to B(H)$ we can define a ``conjugate" representation $\pi^\vee\colon \C[\Irr(\CC)]\to B(\bar H)$ by $\pi^\vee(x)\bar\xi=\overline{\pi(x^\vee)^*\xi}$. For $\pi=\pi_Z$ as above, this is the representation on the Hilbert space $\Mor_{\indC}(Z,\un)$ where $\pi^\vee_Z([X])\xi$ is given by the composition
\begin{multline*}
\bar{R}_X^* (\iota_X \otimes \xi \otimes \iota_{\bar{X}})(c_X^* \otimes \iota_{\bar{X}}) (\iota_Z \otimes \bar{R}_X)\colon Z \xrightarrow{\iota_Z \otimes \bar{R}_X}  Z \otimes X \otimes \bar{X}\xrightarrow{c_X^* \otimes \iota} X \otimes Z \otimes \bar{X} \\
\xrightarrow{\iota\otimes\xi\otimes\iota} X \otimes \bar{X} \xrightarrow{\bar R_X^*} \un.
\end{multline*}
The representation $\pi^\vee_Z$ extends to $C^*(\CC)$, or equivalently, the anti-automorphism $x\mapsto x^\vee$ extends to~$C^*(\CC)$. This is proved in \cite{MR3447719}*{Lemma~6.2} using tube algebras. Another way of seeing this is by using that the class of representations $\pi_Z$ coincides with the class of admissible representations considered in~\cite{MR3406647}. Since admissibility means positivity of certain endomorphisms of $U\otimes\bar U$, by using that the anti-automorphism $T\mapsto T^\vee$ of $\CC(U\otimes\bar U)$ preserves positivity it is easy to check that admissibility is preserved under the operation $\pi\mapsto\pi^\vee$.

\subsection{Ocneanu's tube algebra}

The second algebra associated with $\CC$ is the \emph{tube algebra} $\Tub(\CC)$. As a space,
$$
\Tub(\CC)=\bigoplus_{i,j\in\Irr(\CC)}\Tub(\CC)_{ij},\ \ \Tub(\CC)_{ij}=\bigoplus_{s\in\Irr(\CC)}\CC(U_s\otimes U_j,U_i\otimes U_s).
$$
We denote by $x^s_{ij}$ the component of $x\in\Tub(\CC)$ lying in $\CC(U_s\otimes U_j,U_i\otimes U_s)\subset\Tub(\CC)_{ij}$. Then the product and involution on $\Tub(\CC)$ are defined by
\begin{align*}
(xy)^s_{ij} &=\sum_{\substack{k,r,t\in\Irr(\CC),\\w\in\onb\CC(U_s,U_r\otimes U_t)}}(\iota_i\otimes w^*)(x^r_{ik}\otimes\iota_t)(\iota_r\otimes y^t_{kj})(w\otimes\iota_j),\\
(x^*)^s_{ij} &=(\bar R_s^*\otimes\iota_i\otimes\iota_s)(\iota_s\otimes (x^{\bar s}_{ji})^*\otimes\iota_s)(\iota_s\otimes\iota_j\otimes R_s),
\end{align*}
where $\onb\CC(U_s,U_r\otimes U_t)$ denotes an orthonormal basis in $\CC(U_s,U_r\otimes U_t)$, and we take the dual of $\bar U_s$ among the $(U_t)_t$, denoted by $U_{\bar s}$.

The tube algebra has a $*$-anti-automorphism $x\mapsto x^\vee$ defined by
$$
(x^\vee)^s_{ij}=(x^{\bar s}_{\bar{\jmath}\bar{\imath}})^\vee.
$$
Here, in order to compute $(x^{\bar s}_{\bar{\jmath}\bar{\imath}})^\vee$, we use the solutions of the conjugate equations for the tensor products obtained from those for the factors: $R_{U\otimes V}=(\iota\otimes R_U\otimes \iota)R_V$, $\bar R_{U\otimes V}=(\iota\otimes \bar R_V\otimes \iota)\bar R_U$. Observe that the anti-automorphism $x\mapsto x^\vee$ depends on the choice of the standard solutions $(R_k,\bar R_k)$, but any two such anti-automorphisms differ by a gauge automorphism~$\gamma_z$, $z=(z_i)_i\in\T^{\Irr(\C)}$, defined by $\gamma_z(x)^s_{ij}=z_i\bar z_j x^s_{ij}$. Similarly, the automorphism $x\mapsto x^{\vee\vee}$ is not the identity in general, but a gauge automorphism.

As opposed to the fusion algebra, for any representation of the tube algebra on a pre-Hilbert space the elements of $\Tub(\CC)$ act by bounded operators, with universal bounds on the norms; in particular, the tube algebra admits a universal C$^*$-completion $C^*(\Tub(\CC))$, see \cite{MR3447719}*{Lemma~4.4} or~\cite{arXiv:1511.07329}*{Lemma~3.9}. Similarly to the case of fusion algebras, any object $(Z,c)$ of $\ZC$ defines representations $\pi_Z$ and $\pi^\vee_Z$ of~$\Tub(\CC)$. The second representation is a bit easier to describe. The underlying space is
$$
H^\vee_Z=\bigoplus_{i\in\Irr(\CC)}H^\vee_{Z,i}, \ \ H^\vee_{Z,i}=\Mor_{\indC}(Z,U_i),
$$
and for $x\in\Tub(\CC)$ and $\xi=(\xi_i)_i\in H^\vee_Z$ we have
\begin{equation} \label{eq:pi-vee}
\big(\pi^\vee_Z(x)\xi\big)_i=\sum_{j,s}\sqrt{\frac{d_i}{d_j}} (\iota_i \otimes \bar{R}_s^*)(x^s_{ij}\otimes\iota_{\bar s})(\iota_s\otimes\xi_j\otimes\iota_{\bar s})(c^*_s\otimes\iota_{\bar s})(\iota_Z\otimes \bar R_s).
\end{equation}

Then the representation $\pi_Z$ is defined by letting $H_Z=\overline{H^\vee_Z}$ and $\pi_Z(x)\bar\xi=\overline{\pi^\vee_Z(x^{\vee})^*\xi}$. Expanding the definitions, we get
$$
H_Z=\bigoplus_{i\in\Irr(\CC)}H_{Z,i}, \ \ H_{Z,i}=\Mor_{\indC}(U_i,Z),
$$
and, for $x\in\Tub(\CC)$ and $\xi=(\xi_i)_i\in H_Z$,
\begin{multline*}
\big(\pi_Z(x)\xi\big)_i=\sum_{j,s}\sqrt{\frac{d_i}{d_j}}(\iota_Z\otimes \bar R^*_s)(c_s\otimes\iota_{\bar s}) (\iota_s\otimes\xi_j\otimes\iota_{\bar s})\\ (R^*_{\bar{\imath}}\otimes\iota_s\otimes\iota_j\otimes\iota_{\bar s})(\iota_i\otimes x^s_{\bar{\imath}\bar{\jmath}}\otimes\iota_j\otimes\iota_{\bar s})(\iota_i\otimes\iota_s\otimes\bar R_{\bar{\jmath}}\otimes\iota_{\bar s})(\iota_i \otimes \bar{R}_s).
\end{multline*}
This representation depends on the choice of standard solutions, because the anti-auto\-mor\-phism $x\mapsto x^\vee$ does. But since the gauge automorphisms $\gamma_z$ are unitarily implemented on the space $H_Z$, the equivalence class of~$\pi_Z$ does not depend on any choices.

Any representation of $\Tub(\CC)$ is equivalent to $\pi_Z$ for some $(Z,c)\in\ZC$~\cite{arXiv:1511.07329}*{Proposition~3.14}, so, ignoring the tensor structure, the C$^*$-category $\ZC$ is equivalent to the representation category of $\Tub(\CC)$. The proof of this result is based on a simple relation between the matrix coefficients of $\pi_Z$, or $\pi^\vee_Z$, and the half-braiding $c$. Namely, by~\eqref{eq:pi-vee}, for any $x^s_{ij}\in\CC(U_s\otimes U_j, U_i\otimes U_s)\subset\Tub(\CC)_{ij}$, $\xi_j\in\Mor_\indC(Z,U_j)$ and $\zeta_i\in\Mor_\indC(Z,U_i)$ we have
$$
(\pi^\vee_Z(x)\xi_j,\zeta_i)\iota_i=\sqrt{\frac{d_i}{d_j}} (\iota_i \otimes \bar{R}_s^*)(x^s_{ij}\otimes\iota_{\bar s})(\iota_s\otimes\xi_j\otimes\iota_{\bar s})(c^*_s\otimes\iota_{\bar s})(\iota_Z\otimes \bar R_s)\zeta^*_i.
$$
Applying the trace $\Tr_i$ to both sides, we get
\begin{equation} \label{eq:matr-coef}
(\pi^\vee_Z(x)\xi_j,\zeta_i)=\frac{1}{\sqrt{d_id_j}}\Tr_{U_i\otimes U_s}(x^s_{ij}(\iota\otimes\xi_j)c_s^*(\zeta^*_i\otimes\iota)).
\end{equation}

The tube algebra is equipped with a faithful positive trace $\tau$ defined by
\begin{equation} \label{eq:trace}
\tau(x)=\sum_i\Tr_i(x^e_{ii}),
\end{equation}
where $e\in\Irr(\CC)$ is the index corresponding to the unit object $\un$, see \cite{MR1782145}*{Proposition~3.2}.\footnote{For fusion categories of sectors, the $*$-algebra $\Tub(\CC)$ defined here is isomorphic to the one defined in~\cite{MR1782145}. The isomorphism is given by
$x^s_{ij}\mapsto(is |(x^s_{ij})^*|sj)^*=(j\bar s |(\iota\otimes\iota\otimes\bar R_s^*)(\iota\otimes x^s_{ij}\otimes\iota)(R_s\otimes\iota\otimes\iota)|\bar s i)$.} As a consequence, the canonical map $\Tub(\CC)\to C^*(\Tub(\CC))$ is injective, so we can consider $\Tub(\CC)$ as a subalgebra of $C^*(\Tub(\CC))$.

For every index $i$, consider the projection $p_i\in\Tub(\CC)$ defined by the identity morphism in $\CC(U_i)=\CC(U_e\otimes U_i,U_i\otimes U_e)\subset\Tub(\CC)_{ii}$. This is the unit in $\Tub(\CC)_{ii}$. By mapping~$[U_s]$ onto the identity morphism in $\CC(U_s)=\CC(U_s\otimes U_e,U_e\otimes U_s)\subset\Tub(\CC)_{ee}$, we can identify $\C[\Irr(\CC)]$ with $p_e\Tub(\CC)p_e=\Tub(\CC)_{ee}$. Then the representation $x\mapsto \pi_Z(x)|_{H_{Z,e}}$ of $\Tub(\CC)_{ee}$ is exactly the representation of $\C[\Irr(\CC)]$ defined by $(Z,c)$ that we considered in the previous subsection. Since any representation of $\Tub(\CC)$ is equivalent to $\pi_Z$, it follows that a representation of $\C[\Irr(\CC)]$ extends to~$C^*(\CC)$ if and only if it arises from a representation of~$\Tub(\CC)$ on a larger space. This gives yet another equivalent way of defining $C^*(\CC)$~\cite{MR3447719}.

\bigskip

\section{Tube algebras of representation categories}
\label{sec:tub-alg-rep-cat}

\subsection{Drinfeld double} In this section we take $\CC$ to be $\Rep G$, the representation category of a compact quantum group $G$. Again we mainly follow the conventions of~\cite{MR3204665}, but let us briefly explain the most important ones for the reader's convenience.

A \emph{finite dimensional representation} $U$ of $G$ is an invertible element of $B(H_U)\otimes C(G)$ such that $(\iota\otimes\Delta)(U)=U_{12}U_{13}$. The \emph{tensor product} of two representations $U$ and $V$ is defined by $U\otimes V=U_{13}V_{23}$. The Hopf $*$-algebra of matrix coefficients of finite dimensional representations of $G$ is denoted by $(\C[G],\Delta)$. By $\Rep G$ we mean the category of finite dimensional unitary representations. The corresponding ind-category is the category of all unitary representations of $G$.

The \emph{contragredient} representation to a finite dimensional representation $U$ is defined~by
$$
U^c=(j\otimes\iota)(U^{-1})\in B(H_U^*)\otimes \C[G],
$$
where $j$ is the canonical anti-isomorphism $B(H_U) \to B(H_U^*)$. When $H_U$ is a Hilbert space, we identify the dual space $H_U^*$ with the complex conjugate Hilbert space $\bar H_U$ so that $j$ becomes $j(T)\bar{\xi} = \overline{T^* \xi}$.

For every finite dimensional representation~$U$ of~$G$, we have a representation $\pi_U$ of the dual algebra $\U(G)=\C[G]^*$ on $H_U$ defined by $\pi_U(\omega)=(\iota\otimes\omega)(U)$. Let us also denote the Woronowicz character $f_1\in\U(G)$ by $\rho$.  Then, given a finite dimensional unitary representation $U$ of $G$, the \emph{conjugate} representation of $U$ is defined by
$$
\bar U=\left(j(\pi_U(\rho)^{1/2})\otimes1\right)U^c\left(j(\pi_U(\rho)^{-1/2})\otimes1\right)\in B(\bar H_U)\otimes \C[G].
$$
This is a unitary representation equivalent to $U^c$, and one has $\pi_{\bar U}(\rho)=j(\pi_U(\rho)^{-1})$. We will usually suppress~$\pi_U$ and simply write $\rho\xi$ for $\xi\in H_U$ instead of~$\pi_U(\rho)\xi$.

The representation $\bar U$ is dual to $U$ in $\Rep G$, and as a standard solution of the conjugate equations we usually take
\begin{equation*}\label{eq:conj}
R_U(1)=\sum_x\bar\xi_x\otimes\rho^{-1/2}\xi_x,\ \ \bar R_U(1)=\sum_x\rho^{1/2}\xi_x\otimes\bar\xi_x,
\end{equation*}
where $\{\xi_x\}_x$ is an orthonormal basis in $H_U$. Obviously these morphisms do not depend on the choice of an orthonormal basis, and we have
$d(U)=\dim_q U=\Tr\pi_U(\rho)$. This construction is most natural, but depending on the context it is sometimes necessary to use other realizations of the dual of ${U}$. Namely, when dealing with the tube algebra we need to stick to a chosen family of irreducible representations~$\{U_s\}_s$, so if~$U_s$ happens to be self-dual, we need to use $U_s$ itself as a model of $\bar{U}_s$, and then the corresponding structure morphisms $R_s, \bar{R}_s\colon \C \to H_s \otimes H_s$ do not have the above form.

When we fix representatives $\{U_s\}_s$ of equivalence classes of irreducible representations of $G$, the representations $\pi_s\colon\U(G)\to B(H_s)$ define an isomorphism $\U(G)\cong\prod_sB(H_s)$. We denote by $c_c(\hat G)\subset\U(G)$ the subalgebra corresponding to the algebraic direct sum $\oplus_sB(H_s)$. More invariantly it can be defined as the Fourier transform of $\C[G]$.

Consider the unitary $W=(U_s)_s\in \Mult(c_c(\hat G)\otimes\C[G])$, so that $(\iota\otimes\omega)(W)=\omega$ for all $\omega\in c_c(\hat G)$. Then
$$
(\iota\otimes\Delta)(W)=W_{12}W_{13}\ \ \text{and}\ \ (\Dhat\otimes\iota)(W)=W_{13}W_{23},
$$
where $\Dhat\colon c_c(\hat G)\to \Mult(c_c(\hat G)\otimes c_c(\hat G))$ is the comultiplication dual to the multiplication on $\C[G]$. The unitary~$W$ is nothing else than the multiplicative unitary operator for $G$.

\smallskip

Denote by $\DD(G)$ the Drinfeld, or quantum, double of $G$, and by $\QDG$ the convolution algebra of compactly supported continuous functions on $\DD(G)$. As a space, $\QDG$ is the algebraic tensor product $c_c(\hat G)\otimes\C[G]$. We write $\omega a$ for the element $\omega\otimes a$ of $\QDG$. The product in $\QDG$ is defined using the rule\footnote{This is different from, for example, \cite{MR3238527}, and reflects our conventions for half-braidings.}
$$
\omega(\cdot\, a_{(1)})a_{(2)}=a_{(1)}\omega(a_{(2)}\cdot),
$$
or in other words,
$$
\omega a=a_{(2)}\omega(a_{(3)}\cdot S^{-1}(a_{(1)})).
$$
The involution is defined using the involutions on $\C[G]$ and $c_c(\hat G)$.

The category of (nondegenerate $*$-preserving) representations, or unitary modules, of $\QDG$ can be identified with $\Dcen(\indcat{(\Rep G)})$. The analogous result for finite dimensional Hopf algebras is well-known, but the case of compact quantum groups is not more difficult. Specifically, suppose that a unitary $\QDG$-module $Z$ is given. Then, since $c_c(\hat{G})$ is a subalgebra of $\QDG$, $Z$ can be regarded as an ind-object of $\Rep G$. Moreover, we define a unitary half-braiding $c$ on $Z$ by
\begin{equation}\label{eq:half-braid}
c_U\colon U\otimes Z\to Z\otimes U,\ \ c_U(\xi\otimes\zeta)=W_{21}(\zeta\otimes\xi)\ \ \text{for}\ \ \xi\in H_U,\ \zeta\in Z,
\end{equation}
where we view $W$ as an element of $\Mult(c_c(\hat G)\otimes\QDG)$. The property that the unitaries~$c_U$ are morphisms in $\indcat{(\Rep G)}$ follows from the identity
$$
\Dhat^\op(\omega)W=W\Dhat(\omega),\ \ \text{for all}\ \ \omega\in c_c(\hat G),
$$
in $\Mult(c_c(\hat G)\otimes\QDG)$. This identity, in turn, holds, since by applying $a\otimes\iota$ and using that $(a\otimes\iota)(W)=a$ for all $a\in\C[G]$, we see that it is equivalent to the defining relations in $\QDG$. On the other hand, the half-braiding condition $c_{U\otimes V}=(c_U\otimes\iota)(\iota\otimes c_V)$ follows from $(\Dhat\otimes\iota)(W)=W_{13}W_{23}$.

Conversely, given any unitary half-braiding $c$ on an object $Z$ in $\indcat{(\Rep G)}$, we have $$c_s(\xi\otimes\zeta)=(\tilde U_s)_{21}(\zeta\otimes\xi)$$ for some unitary $\tilde U_s\in B(H_s\otimes Z)$, and by mapping the matrix coefficients of $U_s\in B(H_s)\otimes\C[G]$ into those of~$\tilde U_s$ we get a representation of $\C[G]$ on $Z$. Together with the action of $c_c(\hat G)$ this gives us an $\QDG$-module structure on $Z$.

Now, for any representation of $\QDG$ on a pre-Hilbert space the elements of $\QDG$ act by bounded operators, with universal bounds on the norms, as $\QDG$ is built of matrix algebras and coefficients of unitary matrices over $\C[G]$. In particular, $\QDG$ admits a universal C$^*$-completion $C^*(\DD(G))$. We remark that this also follows from general theory of algebraic quantum groups~\cite{arXiv:funct-an/9704006} and the knowledge of Haar weights~\cite{MR1059324} for quantum doubles, which show that $\QDG$ is the algebra of matrix coefficients of unitary representations of $\coD(G)$.

The map $\QDG\to C^*(\DD(G))$ is injective, as follows, for example, from the following result.

\begin{lemma}\label{lem:cond-exp}
The formula $E_0(a\omega)=\varphi(a)\omega$, where $\varphi$ is the Haar state on $\C[G]$, defines a faithful conditional expectation $E_0\colon\QDG\to c_c(\hat G)$.
\end{lemma}

\bp It follows from the invariance of $\varphi$ and the defining relations in $\QDG$ that $E_0$ also satisfies $E_0(\omega a)=\varphi(a)\omega$ for all $a\in\C[G]$ and $\omega\in c_c(G)$. Thus $E_0$ is a $c_c(\hat G)$-bimodular projection onto~$c_c(\hat G)$. The complete positivity and faithfulness of $E_0$ follow from the corresponding properties of $\varphi$.
\ep

\begin{remark}
Let $(\hat{\phi}, \hat{\psi})$ be left and right Haar weights of $c_0(\hat{G})$.
Then, since the coproduct of $\QDG$ is just $a \omega \mapsto a_{(1)} \omega_{(1)} \otimes a_{(2)} \omega_{(2)}$, the functionals $(\hat{\phi} E_0, \hat{\psi} E_0)$ are left and right Haar weights on $\QDG$. From the above argument it also follows that $\QDG \to C^*_r(\DD(G))$ is injective.
\end{remark}

Consider now the tube algebra
$$
\Tub(\Rep G)=\bigoplus_{i,j}\Tub(\Rep G)_{ij},\ \ \Tub(\Rep G)_{ij}=\bigoplus_s\Hom_G(H_s\otimes H_j,H_i\otimes H_s).
$$
In order to relate it to $\QDG$ consider also the larger $*$-algebra
$$
\Tub(G)=\bigoplus_{i,j}\Tub(G)_{ij},\ \ \Tub(G)_{ij}=\Tub(\Rep G)_{ij}\otimes B(H_{\bar{\jmath}},H_{\bar{\imath}}).
$$
The algebra structure is defined using that on $\Tub(\Rep G)$ and the composition of operators between the spaces~$H_k$. The involution is defined similarly.

\begin{remark}
In fact, $\Tub(G)$ is a particular example of an \emph{annular algebra} considered in~\cite{MR3447719}. In $\Rep G$ we can consider the objects $Y_s = U_s \otimes H_{\bar{s}}$ for $s \in \Irr(G)$, where the part $H_{\bar{s}}$ represents the trivial representation of $G$ on the Hilbert space $H_{\bar{s}}$. Then we have
$$
\bigoplus_s \Mor_{\Rep G}(U_s \otimes Y_j, Y_i \otimes U_s) = \Tub(\Rep G)_{ij} \otimes B(H_{\bar{\jmath}}, H_{\bar{\imath}}),
$$
so $\Tub(G)$ can be identified with $\mathcal{A}\Lambda$ of~\cite{MR3447719}, where $\Lambda = (Y_s)_s$.
\end{remark}

\begin{lemma}[cf.~\cite{MR3447719}*{Proposition~3.5}]\label{lem:tub-corner}
The $*$-algebra $\Tub(\Rep G)$ is isomorphic to a corner in $\Tub(G)$, and this isomorphism extends to the C$^*$-level.
\end{lemma}

\bp
Choose unit vectors $e_i\in H_{\bar{\imath}}$. Let $m_{ij}\colon H_{\bar{\jmath}}\to H_{\bar{\imath}}$ be the rank-one isometry mapping $e_j$ into $e_i$, and consider the projection $f=\sum_ip_i\otimes m_{ii}\in \Mult(\Tub(G))$, where~$p_i$ are the projections introduced at the end of the previous section. Then $\Tub(\Rep G)\cong f\Tub(G)f$, via the isomorphism mapping $x^s_{ij}$ into $x^s_{ij}\otimes m_{ij}$.

In order to see that this isomorphism extends to the C$^*$-level, it suffices to show that the intersection of the positive cone in $\Tub(G)$, defined as the set of finite sums of elements of the form $x^*x$, with $f\Tub(G)f$ coincides with the positive cone in $f\Tub(G)f$. This follows, for example, from fullness of the projection~$f$ in the the following strong sense: there are, possibly infinitely many, elements $x_\alpha\in\Tub(G)f$ such that $\sum_\alpha x_\alpha x_\alpha^*=1$ in $\Mult(\Tub(G))$. By this we mean that for any given $x\in\Tub(G)$ there exists a finite set~$F$ such that the element $x_F=\sum_{\alpha\in F}x_\alpha x_\alpha^*$ has the property $xx_F=x_Fx=x$, while $xx_\alpha x_\alpha^*=x_\alpha x_\alpha^*x=0$ for $\alpha\not\in F$. Such~$x_\alpha$ can be taken to be of the form $p_i\otimes v$, where $v\in B(H_{\bar{\imath}})$ and $v^*v=m_{ii}$. Then, for any $x\in\Tub(G)f$, we see that
$x^*x=\sum_\alpha(x^*x_\alpha)(x_\alpha^* x)$ is positive in $f\Tub(G)f$.
\ep

\begin{theorem}\label{thm:double}
We have an isomorphism of $*$-algebras $\Tub(G)\cong\QDG$.
\end{theorem}

We will need the following simple lemma.

\begin{lemma}\label{lem:equiv-iso}
Suppose that a homomorphism of C$^*$-algebras $\pi\colon A\to B$ induces an equivalence of the representation categories of $A$ and $B$. Then $\pi$ is an isomorphism.
\end{lemma}

\bp
It is clear that $\pi$ is injective, since otherwise any representation $\theta$ of $A$ such that $\theta|_{\ker\pi}$ is nondegenerate could not arise from a representation of $B$.

Now, if $\theta$ is any representation of $B$ on $H$, the endomorphism algebra of $\theta$ is the commutant $\theta(B)'$. The assumption that $\pi$ induces an equivalence implies that $\theta(B)' = \theta(\pi(A))'$, or equivalently, $\theta(B)'' = \theta(\pi(A))''$. Consider the particular case of $\theta$ such that the image of $\theta$ generates the double dual von Neumann algebra~$B^{**}$. Then the induced homomorphism $\pi^{**} \colon A^{**} \to B^{**}$ is surjective. By the Hahn--Banach theorem $\pi$ must have a dense image, and hence it must be surjective in our C$^*$-setting.
\ep

\bp[Proof of Theorem~\ref{thm:double}] Let us write $\A$ for $\Tub(\Rep G)$. Since the representations~$\pi^\vee_Z$ from the previous section are slightly more convenient to work with than $\pi_Z$, we will apply the anti-automorphism $x\mapsto x^\vee$ and prove an equivalent but different result. Namely, consider the map $B(H_{\bar{\jmath}},H_{\bar{\imath}})\to B(\bar H_{\bar{\imath}},\bar H_{\bar{\jmath}})$, $T\mapsto T^\vee$ defined by $T^\vee\bar\xi=\overline{T^*\xi}$. Together with the anti-automorphism  $x\mapsto x^\vee$ of $\A$, which maps $\A_{ij}$ onto $\A_{\bar{\jmath}\bar{\imath}}$, it defines a $*$-anti-isomorphism
of $\Tub(G)$ onto the algebra
$$
\A_G=\bigoplus_{i,j}\A_{ij}\otimes B(\bar H_j,\bar H_i).
$$
Hence, in order to prove the theorem, it suffices to show that $\A_G\cong\QDG^\op$. We will construct such a canonical isomorphism $\pi$. Therefore the isomorphism in the formulation will be defined uniquely up to a gauge automorphism.

In order to define $\pi\colon\A_G\to\QDG^\op$, let us first take a unitary half-braiding~$(Z, c)$, and consider the corresponding representation $\pi^\vee_Z$ of $\A$ on 
$\oplus_j\Hom_G(Z,H_j)$.
We can then define in the obvious way a representation of $\A_G$ on the Hilbert space
$$
H=\bigoplus_j\Hom_G(Z,H_j)\otimes\bar H_j.
$$
This representation can be denoted by $\pi^\vee_Z\otimes\iota$, but we will simply write $T\xi$ instead of $(\pi^\vee\otimes\iota)(T)\xi$.

The computation of matrix coefficients of this representation can be carried out as follows. Choose an orthonormal basis $\{\xi^i_x\}_x$ in $H_i$. Write $m^i_{xy}$ for the corresponding matrix units in $B(H_i)\subset c_c(\hat G)$, and $u^i_{xy}$ for the matrix coefficients of $U_i$. For an element
$$
T\in \Hom_G(H_s\otimes H_j,H_i\otimes H_s)\otimes B(\bar H_j,\bar H_i),
$$
we denote by $T(y'x'z',xyz)$ the matrix coefficients of $T$, viewed as an operator $H_s\otimes H_j\otimes\bar H_j\to H_i\otimes H_s\otimes\bar H_i$, with respect to the bases $\{\xi^s_x\otimes \xi^j_y\otimes\bar\xi^j_z\}_{x,y,z}$ and $\{\xi^i_{y'}\otimes \xi^s_{x'}\otimes\bar\xi^i_{z'}\}_{y',x',z'}$ in $H_s\otimes H_j\otimes\bar H_j$ and $H_i\otimes H_s\otimes\bar H_i$. Let us also write $\rho_s$ for the Woronowicz character $\rho$ acting on $H_s$, and $\rho_{s,x'x}$ for the corresponding matrix coefficients.
Then, by~\eqref{eq:matr-coef}, for $T^s_{ij}\in \CC(U_s\otimes U_j,U_i\otimes U_s)\otimes B(\bar H_j,\bar H_i)$, $\xi_j\in\Hom_G(Z,H_j)$ and $\zeta_i\in\Hom_G(Z,H_i)$, we have
\begin{multline} \label{eq:rep}
(T^s_{ij}(\xi_j\otimes\bar \xi^j_z),\zeta_i\otimes\bar \xi^i_{z'})\\ =\frac{1}{\sqrt{d_id_j}}\sum_{x,y,x',y',x'',y''}\rho_{i,y''y'}\rho_{s,x''x'}T^s_{ij}(y'x'z',xyz)
((\iota\otimes\xi_j)c^*_s(\zeta_i^*\otimes\iota)(\xi^i_{y''}\otimes\xi^s_{x''}),\xi^s_x\otimes\xi^j_y).
\end{multline}

On the other hand, we can also turn $H$ into a unitary $\QDG^\op$-module. Namely, we define the action of $\QDG^\op$ on $\bar Z$ by $X\bar\xi=\overline{X^*\xi}$, and then transfer it to $H$ using the unitary isomorphism
$$
H\cong\bar Z, \ \ \Hom_G(Z,H_j)\otimes\bar H_j\ni \xi\otimes\bar\zeta\mapsto \overline{\xi^*\zeta}.
$$
We remark that then the elements of $c_c(\hat G)^\op$ act in the obvious way,
$$
\omega(\xi\otimes\bar\zeta)=\xi\otimes\overline{\omega^*\zeta}.
$$
The action of $\C[G]^\op$ can be computed using \eqref{eq:half-braid}, which reads as
\begin{equation*}
c_s(\xi^s_y\otimes\zeta)=\sum_xu^s_{xy}\zeta\otimes\xi^s_x.
\end{equation*}
Therefore, for $\xi_j\in\Hom_G(Z,H_j)$ and $\zeta_i\in\Hom_G(Z,H_i)$, we have
\begin{multline*}
(u^{s*}_{x''x}(\xi_j\otimes\bar\xi^j_y),\zeta_i\otimes\bar\xi^i_{y''})=(u^{s*}_{x''x}\overline{\xi_j^*\xi^j_y},\overline{\zeta_i^*\xi^i_{y''}})
=(\zeta_i^*\xi^i_{y''},u^{s}_{x''x}\xi_j^*\xi^j_y)\\
=(\zeta_i^*\xi^i_{y''}\otimes\xi^s_{x''},c_s(\xi^s_x\otimes\xi_j^*\xi^j_y)).
\end{multline*}
Comparing this with \eqref{eq:rep} we get
$$
(T^s_{ij}(\xi_j\otimes\bar \xi^j_z),\zeta_i\otimes\bar \xi^i_{z'})=\frac{1}{\sqrt{d_id_j}}\sum_{\substack{x,y,x',\\y',x'',y''}}\rho_{i,y''y'}\rho_{s,x''x'}T^s_{ij}(y'x'z',xyz)(u^{s*}_{x''x}(\xi_j\otimes\bar\xi^j_y),\zeta_i\otimes\bar\xi^i_{y''}).
$$

We thus see that any element $T^s_{ij}\in\A_G$ acts on $H$ in the same way as the element
\begin{equation} \label{eq:iso-double}
\frac{1}{\sqrt{d_id_j}}\sum_{x,y,z,x',y',z',x'',y''}\rho_{i,y''y'}\rho_{s,x''x'}T^s_{ij}(y'x'z',xyz)\,m^i_{y''z'}\cdot u^{s*}_{x''x}\cdot m^j_{zy}\in\QDG^\op,
\end{equation}
where $\cdot$ denotes the product in $\QDG^\op$. Since the $\QDG^\op$-modules $H$ as above exhaust all unitary modules up to equivalence and the $*$-algebra $\QDG^\op$ admits a faithful representation, it follows that by denoting the element~\eqref{eq:iso-double} by $\pi(T^s_{ij})$ we get a well-defined $*$-homomorphism~$\pi\colon\A_G\to\QDG^\op$. Furthermore, using that any representation of $\A$ is defined by a unitary half-braiding, Lemma~\ref{lem:tub-corner} implies that the $\A_G$-modules $H$ of the above form exhaust all unitary $\A_G$-modules up to equivalence. Therefore the representation categories of $\A_G$ and $\QDG^\op$ are both equivalent to $\Dcen(\indcat{(\Rep G)})$, hence to each other, and $\pi$ implements such an equivalence. By Lemma~\ref{lem:equiv-iso} we conclude that $\pi$ defines an isomorphism of the C$^*$-completions of~$\A_G$ and~$\QDG^\op$.

It remains to show that $\pi$ is bijective even before passing to the C$^*$-completions. Recall that we have a faithful positive trace $\tau$ on $\A$ defined by~\eqref{eq:trace}. Then the formula
$$
F(T)=(\tau\otimes\iota)(T)
$$
defines a faithful completely positive map $F\colon \A_G\to\bigoplus_iB(\bar H_i)\cong c_c(\hat G)^\op$, which is clearly $c_c(\hat G)^\op$-bimodular. On the other hand, by Lemma~\ref{lem:cond-exp} we have a faithful conditional expectation $E_0\colon\QDG^\op\to c_c(\hat G)^\op$ defined by
$E_0(a\cdot\omega)=\varphi(a)\omega$, hence a faithful cp $c_c(\hat G)^\op$-bimodular map
$$
E\colon\QDG^\op\to c_c(\hat G)^\op
$$
defined by $E(a\cdot\omega)=d_i\varphi(a)\omega$ for $a\in\C[G]^\op$ and $\omega\in B(H_i)^\op\subset c_c(\hat G)^\op$.
Then $\pi$ intertwines~$F$ with~$E$. As~$F$ is faithful, this already implies that $\pi$ is injective. Since the image of~$\pi$ is dense in~$\QDG^\op$ in the norm defined by the $c_c(\hat G)^\op$-valued inner product $\langle X, Y\rangle=E(X^*\cdot Y)$, and  the spaces $$B(H_i)\cdot\operatorname{span}\{u^s_{x'x}\}_{x',x}\cdot B(H_j)$$ are mutually orthogonal for different triples $(i,s,j)$. It then follows that the image of $\Hom_G(H_s\otimes H_j,H_i\otimes H_s)\otimes B(\bar H_j,\bar H_i)$ under $\pi$ must coincide with $B(H_i)\cdot\operatorname{span}\{u^{\bar s}_{x'x}\}_{x',x}\cdot B(H_j)$.
\ep

\subsection{Temperley--Lieb categories}

Consider the representation category of $\SU_q(2)$ for $q\in(0,1)$. It is known to be equivalent to the Temperley--Lieb category $\TL(q+q^{-1})$ generated by one object~$U_{1/2}$ and one morphism $R\colon\un\to U_{1/2}\otimes U_{1/2}$ such that
$$
\|R\|^2=q+q^{-1}\ \ \text{and}\ \ (R^*\otimes\iota)(\iota\otimes R)=-\iota,
$$
see, e.g.,~\cite{MR3204665}*{Section~2.5}. The representation theory of the quantum double of $\SU_q(2)$ was studied by Pusz~\cite{MR1213303}, and we can use his results to get information on the structure of $C^*(\Tub(\TL(d)))$ for $d>2$.

The representations of $C^*(\DD(\SU_q(2)))$ are described as follows. As a set, the primitive spectrum is the disjoint union of the closed interval $S_0=[-1,1]$ and countably many circles $S_p=\T$, $p=\frac{1}{2},1,\frac{3}{2},\dots$. The topology is not discussed in~\cite{MR1213303}, but one can at least easily see that on every set $S_p$ we get the standard topology. There are two one-dimensional representations corresponding to the points $\pm1\in S_0$, all other representations are infinite-dimensional. Furthermore, for every $x\in S_p$, excluding the case $x=\pm1$ when $p=0$, the restriction of the representation of $C^*(\DD(\SU_q(2)))$ corresponding to $x$ to $c_0(\widehat{\SU_q(2)})$ decomposes into a direct sum of the spin $i$ modules $H_i$ for $i=p,p+1,\dots$, each appearing with multiplicity one. One immediate consequence of this description is that the C$^*$-algebra $C^*(\DD(\SU_q(2)))$ is liminal, and in particular, of type I.

Consider the closed ideal $J$ in $C^*(\DD(\SU_q(2)))$ generated by the unit $1_0\in B(H_0)=\C$ in the block of $c_0(\widehat{\SU_q(2)})$ corresponding to the spin $0$ (trivial) representation. Then by the above description, the primitive spectrum of $J$ is the interval $S_0$, and the image of~$1_0$ in every irreducible representation of $J$ is a rank one projection. It follows that $J$ is a continuous trace C$^*$-algebra strongly Morita equivalent to $C(S_0)$. In order to describe it explicitly, we have to look more closely at how the representations are defined. This is more transparently done in~\cite{MR1804864} and~\cite{arXiv:1410.6238} using parabolic induction. The underlying algebraic $c_c(\widehat{\SU_q(2)})$-module is $\oplus^\infty_{i=0} H_i$ (only integral spins appear), on which one defines an action of $\C[\SU_q(2)]$ using formulas that only mildly depend on the parameter $x\in S_0$, inducing the structure of a $\QDS$-module. Then one defines a pre-scalar product on $\oplus^\infty_{i=0} H_i$ by rescaling the obvious scalar product on every subspace~$H_i$. This rescaling does not change the scalar product on $H_0\cong\C$, but for all other~$i$ the scalar product  on $H_i$ is multiplied by a factor that tends to $0$ as $x$ converges to $\pm1$, see, for example,~\cite{arXiv:1410.6238} or~\cite{MR2803790}. From this it becomes clear that as an $J$-$C(S_0)$-imprimitivity bimodule we can take the C$^*$-Hilbert $C(S_0)$-module consisting of continuous maps $\xi\colon S_0\to \oplus^\infty_{i=0}H_i$ (where we now consider the Hilbert space direct sum) such that $\xi(\pm1)\in H_0$. It follows that $J$ is isomorphic to the C$^*$-algebra of continuous functions $f$ on $S_0$ with values in the algebra $K(\oplus^\infty_{i=0}H_i)$ of compact operators on $\oplus^\infty_{i=0}H_i$ such that $f(\pm1)\in\C 1_0$.

The description of $C^*(\DD(\SU_q(2)))/J$ is even easier, since for this C$^*$-algebra all irreducible representations are infinite dimensional and, moreover, for all $x\in S_p$, $p=\frac{1}{2},1,\frac{3}{2},\dots$, the underlying space of the representation corresponding to $x$ can be identified with $\oplus^\infty_{i=p}H_i$  (the sum is over $i=p,p+1,\dots$). Specifically, $C^*(\DD(\SU_q(2)))/J$ is isomorphic to the direct sum of the algebras of continuous functions on $S_p$ with values in the compact operators on $\oplus^\infty_{i=p}H_i$.

Therefore the results of Pusz can be summarised as follows.

\begin{proposition}
For every $q\in(0,1)$, we have a short exact sequence
$$
0\to J\to C^*(\DD(\SU_q(2)))\to\bigoplus_{p\in\frac{1}{2}\mathbb N}C(\T)\otimes K\Biggl(\bigoplus^\infty_{i=p}H_i\Biggr)\to0,
$$
where
$$
J=\left\{f\in C\Biggl([-1,1], K\biggl(\bigoplus^\infty_{i=0}H_i\biggr)\Biggr)\mid f(\pm1)\in\C1_0\right\}.
$$
\end{proposition}

By Theorem~\ref{thm:double} and Lemma~\ref{lem:tub-corner}, the tube algebra is obtained by choosing a unit vector in every spin module $H_i$ and then cutting down $C^*(\DD(\SU_q(2)))$ by the sum of the corresponding projections.

\begin{corollary}
For every $d>2$, we have a short exact sequence
$$
0\to {\mathcal J}\to C^*(\Tub(\TL(d)))\to\bigoplus_{i\in\frac{1}{2}\mathbb N}C(\T)\otimes K(\ell^2(\{i,i+1,\dots\}))\to0,
$$
where
$$
{\mathcal J}=\left\{f\in C\left([-1,1], K(\ell^2(\Z_+))\right)\mid f(\pm1)\in\C1_0\right\}.
$$
\end{corollary}

We stress that in this formulation the choice of indices is consistent with the grading on the underlying spaces of the irreducible representations of the tube algebra, so that any element of $\Tub(\TL(d)))_{kl}$ ($k,l\in\frac{1}{2}\Z_+$) maps the $l$-th basis vector in $\ell^2(\{i,i+1,\dots\})$, $\ell^2(\Z_+)$, or $\ell^2(\{0\})$, into a scalar multiple of the $k$-th vector, provided they are both present in the space, and acts as zero otherwise.

\smallskip

The above short exact sequences do not fully describe the hull-kernel topology on the primitive spectrum
$$
\bigcup_{p\in\frac{1}{2}\Z_+}S_p=[-1,1]\sqcup\bigsqcup_{p\in\frac{1}{2}\N}\T
$$
of $C^*(\DD(\SU_q(2)))$, or equivalently, of $C^*(\Tub(\TL(q+q^{-1})))$. But the topology is also not difficult to understand. Since the relative topologies on the subsets $S_0$ and $\cup_{p\in\frac{1}{2}\N}S_p$ are the obvious ones, and the first set is open and the second is closed, the only question is how to describe the closure of $S_0$. The answer is that when points in $S_0\setminus\{z_0^+:=1\}$, resp. $S_0\setminus\{z_0^-:=-1\}$, approach  $z_0^+$, resp. $z_0^-$, they simultaneously approach $z_1^+:=1\in S_1$, resp. $z_1^-:=-1\in S_1$. This corresponds to convergence of complementary series representations (parameterised by points of $S_0\setminus\{z_0^-,z^+_0\}$ close to $z_0^-$ or $z^+_0$) to two particular principal series representations. In other words, we have the following description.

\begin{proposition}
For $q\in(0,1)$, the topology on the primitive spectrum of the C$^*$-algebra $C^*(\DD(\SU_q(2)))$ is described as follows: a set $U$ is open if and only if it is open in the usual topology on the disjoint union of the sets $S_p$ and if the point $1\in S_1$, resp. $-1\in S_1$, lies in $U$, then $U$ must also contain a punctured neighbourhood of $1\in S_0$, resp. of $-1\in S_0$.
\end{proposition}

\bp
In order to prove the proposition we do not need to know how the representations are defined, only how one identifies the set of equivalence classes of irreducible representations that are isomorphic to  $\oplus^\infty_{i=p}H_i$ as representations of~$\SU_q(2)$  with $S_p$. These representations are parameterised by the values of the quantum Casimir~\cite{MR1213303}. For $p=0$ the set of possible values is the interval $[-\sqrt{1+q^2},\sqrt{1+q^2}]$, which we identify with $S_0=[-1,1]$. For $p\in\frac{1}{2}\N$ this set is the ellipse
$$
\left\{z\in\C: \left|\frac{\sqrt{1+q^2}}{q}z-2\right|+\left|\frac{\sqrt{1+q^2}}{q}z+2\right|=2(q^p+q^{-p})\right\},
$$
which we identify with $S_p=\T$ in such a way that the points $\pm\frac{q(q^p+q^{-p})}{\sqrt{1+q^2}}$ get identified with $\pm1$. From this we see that the only pairs of points in $\cup_{p\in\Z_+}S_p$ not separated by the values of the Casimir are $(z^+_0,z^+_1)$ and $(z_0^-,z_1^-)$. Since the closure of $S_0$ must be contained in $\cup_{p\in\Z_+}S_p$ and it is not difficult to see that the Casimir defines a continuous function on the spectrum, this already shows that the only difference from the usual topology on $\cup_{p\in\frac{1}{2}\Z_+}S_p$ that might occur is that a net in $S_0$ converges to $z^+_0$, resp. $z^-_0$, and simultaneously to $z^+_1$, resp. $z^-_1$.

Take a rank one projection $e_1$ in $B(H_1)\subset C^*(\DD(\SU_q(2)))$. We then have a short exact sequence
$$
0\to C_0(S_0\setminus\{z^+_0,z^-_0\})\to e_1C^*(\DD(\SU_q(2)))e_1\to C(S_1)\to0.
$$
It follows that the unital C$^*$-algebra $e_1C^*(\DD(\SU_q(2)))e_1$ is abelian and its spectrum is $(S_0\setminus\{z^+_0,z^-_0\})\cup S_1$. The topology is inherited from the spectrum of $C^*(\DD(\SU_q(2)))$, since $e_1C^*(\DD(\SU_q(2)))e_1$ is strongly Morita equivalent to the ideal of $C^*(\DD(\SU_q(2)))$ generated by $e_1$. Since $(S_0\setminus\{z^+_0,z^-_0\})\cup S_1$ is compact, we see that if a net in $S_0\setminus\{z^+_0,z^-_0\}$ converges to $z^+_0$, resp. $z^-_0$, then in $(S_0\setminus\{z^+_0,z^-_0\})\cup S_1$ is has no other choice than to converge to~$z^+_1$, resp.~$z^-_1$.

Finally, it is clear that the points $z^+_0$ and $z^-_0$ are both closed in the primitive spectrum of $C^*(\DD(\SU_q(2)))$.
\ep

Note that this proposition and its proof show that if for $p\in\frac{1}{2}\Z_+$ we choose a rank one projection $e_p\in B(H_p)\subset C^*(\DD(\SU_q(2)))$, then the C$^*$-algebra $e_pC^*(\DD(\SU_q(2)))e_p$ will be abelian, with spectrum $X_p$, where $X_0=S_0$, $X_p=S_{\frac{1}{2}}\cup S_{\frac{3}{2}}\cup\dots\cup S_{p}$ for $p\in\frac{1}{2}\N$, and $X_p=(S_0\setminus\{z^+_0,z^-_0\})\cup S_1\cup\dots\cup S_p$ for $p\in\N$, where in the last case the topology on $X_p$ is obtained by gluing $S_0$ to $S_1$ via the identification of~$z^\pm_0$ with~$z^\pm_1$. This recovers the description of $\Tub(\TL(q+q^{-1}))_{pp}$ obtained by Ghosh and Jones~\cite{MR3447719}.

\bigskip

\section{Morita equivalence of categories and tube algebras}
\label{sec:mor-eqv-cat-tub-alg}

\subsection{\texorpdfstring{$2$}{2}-categories and \texorpdfstring{$Q$}{Q}-systems}

Let $\BB$ be a small C$^*$-$2$-category. The definition is a straightforward adaptation of the standard algebraic one~\cite{MR1712872}, but as in~\cite{MR1887093}, we denote by $\otimes$ the horizontal composition of morphisms. Thus $\BB$ is given by
\begin{itemize}

\item a set $\Lambda$  ($0$-cells);

\item small C$^*$-categories $\BB_{st}$ for all $s,t\in \Lambda$ (the objects of $\BB_{st}$ are called $1$-morphisms $t\to s$);

\item bilinear unitary bifunctors $\otimes\colon\BB_{rs}\times\BB_{st}\to\BB_{rt}$ and unit objects $\un_s\in\BB_{ss}$.
\end{itemize}
The axioms which this structure should satisfy are analogous to those of strict C$^*$-tensor categories. In other words, the main difference from the latter categories is that the tensor product $X\otimes Y$ is defined not for all objects but only when $X\in \BB_{rs}$ and $Y\in \BB_{st}$, and is in $\BB_{r t}$.

We also always assume that the units $\un_s$ are simple and the categories $\BB_{st}$ are closed under subobjects and finite direct sums.

\smallskip

We are mainly interested in rigid C$^*$-$2$-categories with the set $\Lambda$ of $0$-cells consisting of two points. Such categories can equivalently be described in terms of pairs $(\CC,Q)$ consisting of a rigid C$^*$-tensor category $\CC$ (satisfying our usual assumptions) and a standard simple $Q$-system $Q$ in $\CC$. This means that we are given an isometry $v\colon\un\to Q$ and an isometric up to a scalar factor morphism $w\colon Q\to Q\otimes Q$ such that $(Q,w^*,v)$ is an algebra in $\CC$. Recall that then the Frobenius compatibility condition $(w^*\otimes\iota)(\iota\otimes w)=ww^*$ is satisfied~\citelist{\cite{MR1444286}\cite{MR3308880}}. The assumptions of simplicity and standardness mean that~$Q$ is simple as a $Q$-bimodule and $w^*w=d(Q)\iota$. We will sometimes write $m_Q$ for the product $w^*\colon Q\otimes Q\to Q$. Let us briefly describe this correspondence. 

First, given a pair $(\CC,Q)$ as above, we can construct a C$^*$-$2$-category $\BB^Q$ of unitary modules in $\CC$ in the following way~\citelist{\cite{MR1887093}\cite{MR1966524}}. Take $\Lambda=\{0,1\}$, $\BB^Q_{00}=\CC$, $\BB^Q_{11}=\QmodQ$, $\BB^Q_{01}=\modQ$ and $\BB^Q_{10}=\Qmod$. The tensor products are defined over~$Q$ when possible, otherwise they are taken in $\CC$. To be more precise, depending on the chosen model, the tensor product over $Q$ may not be strictly associative, so we get a bicategory rather than a $2$-category, which then has to be strictified. As is common, we are going to ignore this minor issue. Let us note also in passing that, as in~\cite{arXiv:1501.07390}, we normalise the structure morphisms $P_{M,N}\colon M\otimes N\to M\otimes_Q N$ so that $P_{M,N}P_{M,N}^*=d(Q)\iota$. For any left $Q$-module~$M$ this allows us to take~$M$ as a model of $Q\otimes_QM$, with $P_{Q,M}$ given by the morphism $m^l_M\colon Q\otimes M\to M$ defining the module structure. Of course, the same can be said for the right modules.

The $2$-category $\BB^Q$ is known to be rigid~\citelist{\cite{MR2075605}\cite{arXiv:1501.07390}}. Let us give explicit formulas for standard solutions of the conjugate equations for one-sided modules, which were omitted in~\cite{arXiv:1501.07390}. Since an object of $\CC$ can be considered as an object in different components~$\BB^Q_{st}$ of~$\BB^Q$, let us continue to denote the dimension function on $\BB^Q_{00}=\CC$ by~$d$, and denote the dimension functions on $\Qmod$, $\modQ$ and $\QmodQ$ by $d^{\Qmod}$, $d^{\modQ}$ and $d^Q$. By \cite{arXiv:1501.07390}*{Proposition~6.9} we have $d^Q(M)=d(Q)^{-1}d(M)$. By virtue of multiplicativity and invariance under passing to the dual of the dimension function on $\BB^Q$, this implies that $d^{\Qmod}(M)=d(Q)^{-1/2}d(M)$ and similarly for $d^{\modQ}$. Now, given a left $Q$-module~$M$ (so $M\in\BB^Q_{10}$), we have a solution given in~\cite{arXiv:1501.07390}*{Lemma~6.6}. By rescaling it we get the solution
\begin{align*}
R^{\Qmod}_M = d(Q)^{-3/4} P_{M,\bar{M}} {R}_M &\colon \un \to \bar M \otimes_Q {M},\\
\bar R^{\Qmod}_M = d(Q)^{-1/4} (m^l_M\otimes \iota_{\bar M})(\iota_Q\otimes \bar R_M)&\colon Q \to M \otimes \bar M,
\end{align*}
where $(R_M, \bar{R}_M)$ is a standard solution of the conjugate equations for $M$ in $\CC$. Therefore
$$
(\bar R^{\Qmod *} \otimes_Q \iota_{{M}}) (\iota_{{M}} \otimes {R}^{\Qmod}_M) = \iota_{{M}}, \quad (\iota_{\bar M}\otimes_Q\bar{R}^{\Qmod *}_M) (R^{\Qmod}_M\otimes\iota_{\bar M}) = \iota_{\bar M}.
$$
Using the considerations in~\cite{arXiv:1501.07390}*{Section~6.2} it is not difficult to check that $$\norm{R^{\Qmod}_M}^2=\norm{\bar{R}^{\Qmod}_M}^2 = d(Q)^{-1/2} d(M) = d^{\Qmod}(M),$$ so the solution $(R^{\Qmod}_M,\bar R^{\Qmod}_M)$ is standard.

Conversely, given a rigid C$^*$-$2$-category $\BB$ with two $0$-cells,  we can construct a $Q$-system in $\CC=\BB_{00}$ by taking any nonzero object $X\in\BB_{10}$ and letting $Q=\bar X\otimes X$, $v=d(X)^{-1/2}R_X$ and $w=d(X)^{1/2}\iota_{\bar X}\otimes\bar R_X\otimes\iota_X$. By a result of M{\"u}ger~\cite{MR1966524}*{Proposition~4.5}, the $2$-categories $\BB$ and $\BB^Q$ are equivalent. Explicitly, the equivalence is given by the functor
$$
\BB\to\BB^Q,\ \ \begin{pmatrix} Y_{00} & Y_{01}\\ Y_{10}\ & Y_{11}\end{pmatrix}\mapsto \begin{pmatrix} Y_{00} & Y_{01}\otimes X\\ \bar X\otimes Y_{10}\ & \bar X\otimes Y_{11}\otimes X\end{pmatrix}.
$$
Here, for example, the left $Q$-module structure on $\bar X\otimes Y$ is given by
$$
m^l_{\bar X\otimes Y}=d(X)^{1/2}\iota_{\bar X}\otimes\bar R_X^*\otimes\iota_X\colon (\bar X\otimes X)\otimes(\bar X\otimes Y)\to\bar X\otimes Y.
$$
Furthermore, if we took a simple object $X$, then the $Q$-system $Q=\bar X\otimes X$ would satisfy $\dim\CC(\un,Q)=1$, which is equivalent to \emph{irreducibility} of $Q$, meaning that $Q$ would be simple as a left and right $Q$-module.

\smallskip

Following M{\"u}ger~\cite{MR1966524} we say that two rigid C$^*$-tensor categories $\CC$ and $\DD$ are \emph{weakly (unitarily) monoidally Morita equivalent} if there exists a rigid C$^*$-$2$-category $\BB$ with two $0$-cells such that $\BB_{00}\cong \CC$, $\BB_{11}\cong\DD$ and $\BB_{01}\ne0$ (the rigidity then implies that also $\BB_{10} \ne 0$). By the above discussion this is equivalent to the existence of a standard simple $Q$-system $Q\in\CC$ such that $\DD\cong\QmodQ$. Furthermore, then $Q$ can be chosen to be irreducible.

\begin{remark}
Historically, rigid C$^*$-$2$-categories first appeared in the study of subfactors. Namely, let $N \subset M$ be an extremal II$_1$ subfactor of finite index. Then one looks at the C$^*$-$2$-category $\BB$ of Hilbert $M$-bimodules, $N$-bimodules, $M$-$N$-modules, and $N$-$M$-modules generated by $X = {}_NL^2(M)_M\in\BB_{01} \subset N\mhyph\modl\mhyph M$. Then $Q = {}_NL^2(M)_N$ is the $Q$-system in $\BB_{00} \subset N\mhyph\mathrm{bimod}$ corresponding to the generator $X$, and the product structure of $Q$ comes from that of $M$. Any rigid C$^*$-$2$-category $\BB=\BB^Q$ defined by an irreducible $Q$-system $Q\in\CC$ such that $\CC$ is generated by $Q$, appears in this way, up to unitary equivalence. However, as in~\citelist{\cite{MR1822847}\cite{MR1966524}}, in what follows we do not require this generation property. For example, in the case $Q = \un_\CC$ we obtain a $2$-category with $\BB_{st} = \CC$ for all~$s$,~$t$.
\end{remark}

\subsection{Tube algebras and imprimitivity modules}

Let $\BB$ be a small rigid C$^*$-$2$-category and $\Lambda$ be its set of $0$-cells. For all $s,t\in\Lambda$ choose representatives $U_{sti}$, $i\in I_{st}$, of isomorphisms classes of simple objects in~$\BB_{st}$. We write~$U_{si}$ instead of $U_{ssi}$ and omit the indices $s,t$ altogether when there is only one meaningful choice. We can then define the tube algebra of $\BB$ by
$$
\Tub(\BB)=\bigoplus_{s,t,i,j,k}\BB_{st}(U_k\otimes U_{tj},U_{si}\otimes U_k).
$$
The $*$-algebra structure is defined by exactly the same formulas as for C$^*$-tensor categories. Denote by $\Tub(\BB)^{st}$ the part of $\Tub(\BB)$ corresponding to fixed indices $s,t\in\Lambda$. Then $\Tub(\BB)^{ss}=\Tub(\BB_{ss})$.

Similarly to the case of tensor categories, the representation theory of $\Tub(\BB)$ can be described in terms of the Drinfeld centre $\Dcen(\indcat{\BB})$ of $\indcat{\BB}$, which is defined as follows. An object $(Z,c)$ of $\Dcen(\indcat{\BB})$ is given~by

\begin{itemize}

\item a set $Z=\{Z_s\}_{s\in\Lambda}$ of ind-objects $Z_s\in\indcat{\BB_{ss}}$,

\item collections $c^{st}=(c^{st}_X)_{X\in\BB_{st}}$ of natural unitary isomorphisms $X\otimes Z_t\cong Z_s\otimes X$ satisfying the half-braiding condition
$$
c^{rt}_{X\otimes Y}=(c^{rs}_X\otimes\iota_Y)(\iota_X\otimes c^{st}_Y).
$$
\end{itemize}
This construction seems to be well-known to the experts, and it is discussed in detail in~\cite{MR3398725}.
The centre  $\Dcen(\indcat{\BB})$ is a C$^*$-tensor category in the obvious way.

Every object $(Z,c)\in\Dcen(\indcat{\BB})$ defines representations $\pi_Z$ and $\pi^\vee_Z$ of $\Tub(\BB)$, similarly to how half-braidings define representations of the fusion and tube algebras of C$^*$-tensor categories. Let us write down the full definition only for $\pi^\vee_Z$. The underlying Hilbert space~is
$$
H^\vee_Z=\bigoplus_{s,i}H^\vee_{Z,si},\ \ H^\vee_{si}=\Mor_{\indcat{\BB_{ss}}}(Z_s,U_{si}).
$$
Every element $T\in \Mor(U_k\otimes U_{tj},U_{si}\otimes U_k)\subset\Tub(\BB)$ kills the spaces $H^\vee_{Z,rl}$ with $(r,l)\ne(t,j)$ and maps $\xi\in H^\vee_{Z,tj}$ into the vector in $H^\vee_{Z,si}$ defined as the composition
\begin{multline*}
\sqrt{\frac{d_{si}}{d_{tj}}} (\iota \otimes \bar{R}_{k}^*)(T\otimes\iota)(\iota\otimes\xi\otimes\iota)(c^{st*}_{k}\otimes\iota)(\iota\otimes \bar R_{k})\colon\\
Z_s\to Z_s\otimes U_k\otimes\bar U_k\to U_k\otimes Z_t\otimes\bar U_k\to U_k\otimes U_{tj}\otimes\bar U_k\to U_{si}\otimes U_k\otimes\bar U_k\to U_{si}.
\end{multline*}

\begin{proposition}
The map $(Z,c)\mapsto\pi_Z$ defines a unitary equivalence between the Drinfeld centre $\Dcen(\indcat{\BB})$ and the category of (nondegenerate $*$-preserving) representations of $\Tub(\BB)$.
\end{proposition}

This is proved in exactly the same way as the analogous result for C$^*$-tensor categories~\cite{arXiv:1511.07329}*{Proposition~3.14}, which corresponds to the case when $\Lambda$ consists of one point. One immediate consequence of this proposition is that $\Tub(\BB)$ admits a universal C$^*$-completion $C^*(\Tub(\BB))$.

Denote by $C^*(\Tub(\BB))^{st}$ the closure of $\Tub(\BB)^{st}$ in $C^*(\Tub(\BB))$. Then the space $C^*(\Tub(\BB))^{st}$ is a C$^*$-Hilbert $C^*(\Tub(\BB))^{ss}$-$C^*(\Tub(\BB))^{tt}$-bimodule. If the corners $\BB_{st}$ are nontrivial, then the following theorem shows that $C^*(\Tub(\BB))^{st}$ are $C^*(\Tub(\BB_{ss}))$-$C^*(\Tub(\BB_{tt}))$-imprimitivity bimodules, and the result is true even purely algebraically.

\begin{theorem}\label{thm:alg-imprimitivity}
Assume $\CC$ and $\DD$ are weakly monoidally Morita equivalent rigid C$^*$-tensor categories, and~$\BB$ is the corresponding rigid $C^*$-$2$-category.
Then $\Tub(\BB)^{01}$  is an algebraic imprimitivity $\Tub(\CC)$-$\Tub(\DD)$-bimodule, that is, we have $\Tub(\BB)^{s t} \Tub(\BB)^{t s} = \Tub(\BB)^{s s}$ for  all $s, t \in \{0, 1\}$, and the closure $C^*(\Tub(\BB))^{01}$ of $\Tub(\BB)^{01}$ in $C^*(\Tub(\BB))$ defines a $C^*(\Tub(\BB_{00}))$-$C^*(\Tub(\BB_{11}))$-imprimitivity bimodule.
\end{theorem}

\bp It is enough to prove the first assertion for $t = 0$ and $s = 1$. As discussed in the previous subsection, we may assume that $\BB=\BB^Q$ for an irreducible $Q$-system $Q\in\CC$.

In the proof we will use the following simple observation. Let $M$ be a right $Q$-module and $Y$ be a $Q$-bimodule. Then $P_{M,Y}(m^r_M\otimes\iota)=P_{M,Y}(\iota\otimes m^l_Y)$ by the definition of the tensor product over $Q$. This can be interpreted by saying that if we take $M \otimes Y$ as a model of $(M \otimes Q) \otimes_Q Y$ with the structure morphism $P_{M\otimes Q,Y}=\iota_M \otimes m_Y^l$, then the morphism $m_M^{r} \otimes_Q \iota_Y\colon (M\otimes Q)\otimes_QY=M \otimes Y \to M \otimes_Q Y$ is equal to $P_{M,Y}$. In particular, we have $m_Y^{l} = m_Q\otimes_Q \iota_Y$.

Now, let $\{U_i\}_{i\in I_{00}}$, $\{M_k\}_{k\in I_{01}}$ and $\{X_a\}_{a\in I_{11}}$ be representatives of the isomorphism classes of simple objects, respectively, of $\CC = \BB_{00}$, of the right $Q$-modules $\modQ = \BB_{01}$, and of the $Q$-bimodules $\QmodQ = \BB_{11}$, with distinguished choices $U_e = \un$ and $M_e = X_e = Q$. Then we have
\begin{align*}
\Tub(\BB)^{0 1} &= \bigoplus_{i, a, k} \Mor_{\modQ}(M_k \otimes_Q X_a,  U_i \otimes M_k), \\
\Tub(\BB)^{1 0} &= \bigoplus_{i, a, k} \Mor_{\Qmod}(\bar{M}_k \otimes U_i, X_a \otimes_Q \bar{M}_k).
\end{align*}
For each $a \in I_{11}$, denote by $p(a) \in \Tub(\BB)^{11}=\Tub(\BB_{11})$ the unit in $\Tub(\BB_{11})_{aa}$, so
$$
p(a)^b_{c d} = \delta_{e,b} \delta_{a,c}\delta_{a,d} \iota \in \Mor_{\QmodQ}(X_b \otimes_Q X_d, X_c \otimes_Q X_b).
$$
The family $(p(a))_a$ generates $\Tub(\BB)^{11}$ as a two-sided ideal. Therefore it is enough to show that $p = p(a)$ belongs to $\Tub(\BB)^{1 0} \Tub(\BB)^{0 1}$ for any fixed $a$.

Let $\{u_j^\alpha \colon U_j \to X_a\}_{j,\alpha}$ be a complete orthonormal system of isometries, so that we have $\sum_{j,\alpha} u_j^\alpha u_j^{\alpha *} = \iota_a$. We define an element $T_j^\alpha \in \Tub(\BB)^{1 0}$ by
$$
(T_j^\alpha)^k_{bi} = \delta_{e,k} \delta_{j,i} \delta_{a,b} m_a^l (\iota_Q \otimes u_j^\alpha) \in \Mor_{\Qmod}(\bar{M}_k \otimes U_i, X_b \otimes_Q \bar{M}_k),
$$
where we take $X_a$ itself as a model of $X_a \otimes_Q \bar{M}_e$. We claim that $$d(Q) p = \sum_{j,\alpha} T_j^\alpha T_j^{\alpha*}.$$

The only nontrivial component of $T_j^{\alpha *}\in\Tub(\BB)^{0 1}$ is in $\Mor_{\modQ}(M_e \otimes_Q X_a, U_j \otimes M_e)$. It equals
$$
(T_j^{\alpha *})^e_{ja}=(\bar R^{\modQ*}_Q\otimes\iota\otimes\iota)(\iota\otimes_Q (T_j^{\alpha})^{e*}_{aj}\otimes\iota)(\iota\otimes_Q\iota\otimes_Q R^{\modQ}_Q),
$$
where
$$
R^{\modQ}_Q=d(Q)^{-1/4}w\colon Q\to Q\otimes Q,\ \ \bar R^{\modQ}_Q=d(Q)^{1/4}v\colon \un\to Q\otimes_QQ=Q.
$$
Therefore
$$
(T_j^{\alpha *})^e_{ja}=(v^*\otimes\iota\otimes\iota)(\iota\otimes u^{\alpha*}_j\otimes\iota)(m^{l*}_a\otimes\iota)m^{r*}_a=(u^{\alpha*}_j\otimes\iota)m^{r*}_a\colon X_a\to U_i\otimes Q,
$$
where we used that the morphism
$$
\iota\otimes_Qw\colon X_a\otimes_QQ=X_a\to X_a\otimes_Q (Q\otimes Q)=X_a\otimes Q
$$
coincides with $m_a^{r*}$ by the observation at the beginning of the proof.

By definition, the product $T_j^\alpha T_j^{\alpha *}$ in $\Tub(\BB)$ equals
\begin{multline*}
\sum_{b, \beta} (\iota_a \otimes_Q w_b^{\beta *}) (m_a^l \otimes \iota_Q) (\iota_Q \otimes u_j^\alpha u_j^{\alpha *} \otimes \iota_Q) (\iota_Q \otimes m_a^{r*}) (w_b^\beta \otimes_Q \iota_a) \\
\in \bigoplus_b \Mor_{\QmodQ}(X_b \otimes_Q X_a, X_a \otimes_Q X_b),
\end{multline*}
where $\{w_b^\beta\}_{\beta}$ is an orthonormal basis in $\Mor_{\QmodQ}(X_b, Q\otimes Q)$ and we use the identifications $Q \otimes_Q X_a = X_a = X_a \otimes_Q Q$. Taking the sum over $j$ and $\alpha$, we get
$$
\sum_{j,\alpha} T_j^\alpha T_j^{\alpha*} = \sum_{b, \beta} (\iota_a \otimes_Q w_b^{\beta *}) (m_a^l\otimes \iota_Q)  (\iota_Q \otimes m_a^{r*}) (w_b^\beta \otimes_Q \iota_a).
$$
Using again that $m^l_a=m_Q\otimes_Q\iota_a$, this can be written as
$$
\sum_{b,\beta} (w^* w_b^\beta) \otimes_Q ((\iota_a \otimes_Q w_b^{\beta*}) m_a^{r*}).
$$
Since $d(Q)^{-1/2} w$ can be taken as a part of $\{w_b^\beta\}_{b,\beta}$, only the term corresponding $b = e$ with a unique $\beta$ survives, and we obtain $\delta_{e,b} d(Q) \iota_a$ as a result. This proves the claim.

\smallskip

For the last statement of the theorem we need to check that the identification of $\Tub(\BB_{ss})$ with $\Tub(\BB)^{ss}$ extends to an isomorphism of $C^*(\Tub(\BB_{ss}))$ onto the corner $C^*(\Tub(\BB))^{ss}$ of $C^*(\Tub(\BB))$. This is shown in the same way as in the proof of Lemma~\ref{lem:tub-corner}.
\ep

\begin{corollary}
The forgetful functor $\Dcen(\indcat{\BB})\to \Dcen(\indcat{\BB_{ss}})$ is a unitary monoidal equivalence for $s=0,1$.
\end{corollary}

\bp By the theorem, $C^*(\Tub(\BB))$ is the linking algebra between $C^*(\Tub(\BB_{00}))$ and $C^*(\Tub(\BB_{11}))$, in particular, it is itself strongly Morita equivalent to either of these C$^*$-algebras. Therefore the representation categories of $C^*(\Tub(\BB))$ and $C^*(\Tub(\BB_{ss}))$ are equivalent.  Recalling that they are also equivalent to, resp., $\Dcen(\indcat{\BB})$ and $\Dcen(\indcat{\BB_{ss}})$, we conclude that the forgetful functor $\Dcen(\indcat{\BB})\to \Dcen(\indcat{\BB_{ss}})$ is an equivalence of categories. Since this functor is monoidal, this proves the assertion.
\ep

\begin{remark}\label{rem:compar-with-Sch}
This corollary can also be proved in a more direct way as follows. As above, we assume $\BB=\BB^Q$. It is known that $\Dcen(\indcat{\BB_{00}}) = \ZC$ and $\Dcen(\indcat{\BB_{11}}) = \Dcen(\indcat{\QmodQ})$ are monoidally equivalent through Schauenburg's induction $(Z, c) \mapsto (Z \otimes Q, \tilde{c})$, where the left $Q$-module structure on $Z \otimes Q$ is defined by that on $Q\otimes Z$ using the isomorphism $c_Q\colon Q\otimes Z\to Z\otimes Q$, and $\tilde{c}_Y$ is given by $c_Y$ up to the identifications of the form $(Z \otimes Q) \otimes_Q Y = Z \otimes Y$, see~\citelist{\cite{MR1822847}\cite{arXiv:1501.07390}*{Section~6.3}}. The centre of $\BB$ just packages this monoidal equivalence of the centres in the following way.

In order to avoid confusion, let us put $M = Q$ as an object in $\BB_{01} = \modQ$, and $\bar{M} = Q$ as an object in $\BB_{1 0} = \Qmod$. Assume that $(Z,c)$ is an object in the centre of $\indcat{\BB}$. First it gives us $Q$-module isomorphisms $\bar{M} \otimes Z_0 \to Z_1 \otimes_Q \bar{M}$ and $M \otimes_Q Z_1 \to Z_0 \otimes M$. Let $Y$ be a $Q$-bimodule, and put $X = M \otimes_Q Y$ (that is, $Y$ is viewed as an object in $\modQ$). Then we have a commutative diagram
$$
\begin{tikzcd}[column sep=5em]
& & Z_0 \otimes X \otimes_Q \bar{M} \\
X \otimes_Q \bar{M} \otimes Z_0 \ar[r, "\iota_X \otimes_Q c^{1 0}_{\bar{M}}"] \ar[rru, bend left=8, "c^{00}_{X \otimes_Q M}"] \ar[rrd, "\iota_M \otimes_Q c^{10}_{Y\otimes_Q\bar{M}}"'] & X \otimes_Q Z_1 \otimes_Q \bar{M} \ar[ru, "c^{01}_X \otimes_Q \iota_{\bar{M}}"'] \ar[rd, "\iota_M \otimes_Q c^{11}_Y \otimes_Q \iota_{\bar{M}}"] \\
& & M \otimes_Q Z_1 \otimes_Q Y \otimes_Q \bar{M} \ar[uu, "c^{01}_M \otimes_Q \iota_{Y \otimes_Q \bar{M}}"'].
\end{tikzcd}
$$
Taking $Y = Q$, $Y = Q \otimes Q$ and using the naturality for $m_Q$, we see that the $Q$-bimodule structure on~$Z_1$ is identified with the one on $Z_0 \otimes Q$ (see the remark at the beginning of the proof of Theorem~\ref{thm:alg-imprimitivity}). Next, taking~$Y$ to be generic, the commutativity of the upper-left and right triangles above shows that $c^{11}$ is determined by $c^{00}$ (and is given by Schauenburg's induction of $(Z_0,c^{00})$). By considering free one-sided modules $Q\otimes X$ and $X\otimes Q$, we also see that~$c^{01}$ and~$c^{10}$ are completely determined by $c^{00}$, $c^{11}$, $c^{01}_M$, and $c^{10}_{\bar M}$. Thus, up to an isomorphism, $(Z,c)$ is completely determined by $(Z_0,c^{00})$. The argument also shows how given any $(Z_0,c^{00})\in\Dcen(\indcat{\BB_{00}})$ we can construct an object of $\Dcen(\indcat{\BB})$.
\end{remark}

\begin{remark}
Analogues of Theorem~\ref{thm:alg-imprimitivity} for $\C[\CC]$ and $C^*(\CC)$ do not hold in general, already for fusion categories. By~\cite{MR1966525}*{Theorem~4.14}, when $G$ is a finite group, the monoidal category of $G$-bimodules $G\mhyph\mathrm{bimod} = \Rep (G \times G^{\op})$ with usual tensor product over $\C$, is weakly Morita equivalent to $(C(G)\rtimes_{\operatorname{Ad}}G)\mhyph\mathrm{mod}$. Here $C(G) \rtimes_{\operatorname{Ad}}G$ is the crossed product of the function algebra $C(G)$ by the adjoint action of $G$, and is isomorphic to $C^*(\DD(G))$. Specifically, one may take $Q = C(G)$ in $G\mhyph\mathrm{bimod}$, and the equivalence
$$
C(G)\mhyph\mathrm{bimod}_{G\mhyph\mathrm{bimod}} \to (C(G) \rtimes_{\operatorname{Ad}} G) \mhyph\mathrm{mod}
$$
is given by $E \mapsto \oplus_g \delta_g E \delta_e$. Since both $G\mhyph\mathrm{bimod}$ and $(C(G) \rtimes_{\operatorname{Ad}} G) \mhyph\mathrm{mod}$ have commutative fusion rules, their fusion algebras are strongly Morita equivalent if and only if they have the same number of irreducible classes. When $G = S_3$, there are $9$ irreducible classes in $G\mhyph\mathrm{bimod}$, while $C(G)\rtimes_{\operatorname{Ad}}G\mhyph\mathrm{mod}$ has only $8$ irreducible classes, as seen from the computation of conjugacy classes and associated stabiliser subgroups.
\end{remark}

\begin{remark}
Theorem~\ref{thm:alg-imprimitivity} can be used to give a short and transparent proof of \cite{arXiv:1501.07390}*{Corollary~6.19} on invariance of property (T) under weak Morita equivalence, by observing that the representation of $\Tub(\CC)$ that one gets by inducing the representation~$\pi^\vee_{\un_{\Dcen(\DD)}}$ of~$\Tub(\DD)$ is~$\pi^\vee_{\un_{\Dcen(\CC)}}$.
\end{remark}

We finish by constructing a class of half-braidings generalising regular half-braidings that we studied in~\cite{arXiv:1501.07390}. Let $Q$ be  a standard simple $Q$-system in a C$^*$-tensor category $\CC$, and $\BB=\BB^Q$. Let $\{U_i\}_{i\in I_{00}}$, $\{M_k\}_{k\in I_{01}}$ and $\{X_a\}_{a\in I_{11}}$ be representatives of the isomorphism classes of simple objects in $\BB_{00}$, $\BB_{01}$ and~$\BB_{11}$.

\begin{proposition}\label{prop:Zreg-modQ-Y-has-h-br}
Let $Y$ be a $Q$-bimodule. Then there is an object $(Z, c)$ of $\Dcen(\indcat{\BB})$ such that
\begin{align*}
Z_0 = \Zreg(\modQ; Y) &= \bigoplus_k M_k \otimes_Q Y \otimes_Q \bar{M}_k,\\
Z_1 = \Zreg(\QmodQ; Y) &= \bigoplus_a X_a \otimes_Q Y \otimes_Q \bar{X}_a.
\end{align*}
\end{proposition}

\begin{remark}
When $Y = Q$, then $c^{11}$ coincides with the regular half-braiding on $\Zreg(\QmodQ)$ defined in~\cite{arXiv:1501.07390}*{Section~3.2}, while
$(\Zreg(\modQ; Q),c^{00})$ coincides with the object $(\Zreg(\modQ),c)$ which appeared in \cite{arXiv:1501.07390}*{Section~6.5}. For $Q = \un$, analogous constructions carried out through the notion of coends can be found in the literature on fusion categories, see for example~\cite{arXiv:1504.01178}*{Section~3}.
\end{remark}

\bp[Proof of Proposition~\ref{prop:Zreg-modQ-Y-has-h-br}]
Let $X$ be a right $Q$-module, and fix a standard solution of the conjugate equations $(R_X, \bar{R}_X)$ for~$X$ in~$\CC$. Consider the morphism $c_{X, k a}$ from $X \otimes_Q X_a \otimes_Q Y \otimes_Q \bar{X}_a$ to $M_k \otimes_Q Y \otimes_Q \bar{M}_k \otimes X$, given by
$$
c^{01}_{X,k a} = \sqrt{\frac{d^Q_a}{d^{\modQ}_k}}\sum_\alpha (u_{a k}^{\alpha *} \otimes _Q \iota_Y \otimes_Q u_{a k}^{\alpha \vee} \otimes \iota_X) (\iota_{X \otimes_Q X_a \otimes_Q Y \otimes_Q \bar{X}_a} \otimes_Q R^{\modQ}_X),
$$
where $\{u_{a k}^\alpha\}_\alpha$ is an orthonormal basis in $\Mor_{\modQ}(M_k, X \otimes_Q X_a)$. Similarly to~\cite{arXiv:1501.07390}*{Section~3.2}, these morphisms do not depend on the choice of a standard solution for $X$, but they do depend on such a choice for $M_k$ and $X_a$. Collecting them together, we obtain a morphism
$$
c^{01}_X = (c^{01}_{X,k a})_{a, k} \colon X \otimes_Q \Zreg(\QmodQ;Y) \to \Zreg(\modQ;Y) \otimes X.
$$
Then the same argument as in the proof of~\cite{arXiv:1501.07390}*{Lemma~6.21} shows that $c^{01}_X$ is unitary.

The rest of the components of $c$ can be defined in a similar way. For example, $c^{00}$ is defined using standard solutions in $\CC$, while $c^{11}$ is defined with ones in $\QmodQ$. It remains to show that the family $(c_X)_X$ satisfies the half-braiding condition. Again, this can be done in the same way as in the proof of~\cite{arXiv:1501.07390}*{Proposition~6.22}, by reducing the argument to the multiplicativity of standard solutions.
\ep

The objects we have just defined allow us to look from a different angle at some of the previous results and constructions. Namely, using the insight of M{\"u}ger in~\cite{MR1966525}, we may view $\End_{\ZC}(\oplus_i \Zreg(\CC; U_i))$ as a von Neumann algebraic analogue of the tube algebra of $\CC$. Then, in view of Remark~\ref{rem:compar-with-Sch},
$$
\End_{\ZC}\left(\bigoplus_a \Zreg(\modQ; X_a)\right)
$$
corresponds to the tube algebra of $\QmodQ$, while
$$
\Mor_{\ZC}\left(\bigoplus_{i} \Zreg(\CC; U_i), \bigoplus_{a} \Zreg(\modQ; X_a)\right)
$$
gives an imprimitivity bimodule (in the von Neumann algebraic sense) between the two.

\bigskip

\end{document}